\documentclass[12pt,a4paper,titlepage]{article}
\usepackage[latin9]{inputenc}
\usepackage{mathrsfs}
\usepackage{amsmath}
\usepackage[natbibapa]{apacite}

\usepackage{here}
\usepackage{mathrsfs}

\makeatother

\begin{document}

\title{\textbf{Model Averaging for Generalized Linear Model with Covariates that are Missing completely at Random}}

\author{Miaomiao Zheng\thanks{School of Mathematical Sciences, Capital Normal University, Beijing
China. E-mail: 2150502075@cnu.edu.cn.}\textsuperscript{} and Qingfeng Liu\thanks{Corresponding author. Department of Economics, Otaru University of
Commerce, Otaru City, Hokkaido, Japan. E-mail: qliu@res.otaru-uc.ac.jp.
Liu acknowledges financial support from JSPS KAKENHI Grant Number
JP16K03590.}}
\maketitle
\begin{abstract}
In this paper, we consider the estimation of generalized linear models
with covariates that are missing completely at random. We propose
a model averaging estimation method and prove that the corresponding
model averaging estimator is asymptotically optimal under certain
assumptions. Simulaiton results illustrate that this method has better
performance than other alternatives under most situations.

Key Words: Missing data completely at random, generalized linear models,
asymptotically optimal. 
\end{abstract}
%%%%%%%%%%%%%%%%%%%%%%%%%%%%%%%%%%%%%%%%%%%%%%%%%%%%%%%%%%%%%%%%%%%%%%%%%%%%%%%%%%%%%%%%%%%%%%%%%%%%%%%%%%%%%%%%%%%%%%%%%%%%%%%%%%%%%%%%%%%%%%%%%%%%%%%%%%%%%%%%%%%%%%%%%%%%%%%%%%%%%

\section{Introduction}

Our attempt in this paper is at developing an optimal model averaging
method for generalized linear models (GLMs) with missing values on
some covariates. Model averaging, as an alternative of model selection,
is widely adopted for dealing with the uncertainty that the model
selection takes. Unlike Bayesian model averaging, the current paper
focuses on the method of determining weights from frequencies perspective.

The problem of missing values on some covariates is one of the most
common challenge facing empirical researchers, which will bring about
a great impact on the subsequent modeling as well as inference process.
There is a large collection of literature dealing with missing covariates.
See \cite{little1992regression} and \cite{toutenburg2002parametric}
for reviews of this topic. There are three types of missing data:
missing completely at random (MCAR), missing at random (MAR) and missing
not at random. Referring to \cite{little1992regression}, data is
MCAR if the probability of missing or not for a covariate is independent
with any covariates (including itself). Data is MAR if the probability
of missing or not for a covariate is independent with any covariates'
missing values and may depend on the covariates' observed values only.
Except for MCAR and MAR, the other types of missing data are all missing
not at random. In the current paper, we focus on the case that the
covariates are MCAR, which can happen in practice.

A very straightforward and widely used method to handle with missing
covariates problem is complete-case (CC) analysis. The CC method estimates
a model using observations with complete data and discards other observations
that have some missing values of some covariates. Obviously, the discarding
of data seems an unnecessary waste of information. Another method
that is widely used to deal with missing data is mean imputation.
However, it is well known that data imputation generally depends on
observed covariates, so it is not appropriate for the MCAR situation.
In order to make full use of the sample information and avoid the
uncertain influence imputation bring, we propose a model averaging
method, which is based on model averaging for GLMs \cite{zhang2016optimal}
to combine some submodels, including one CC model and some partial
data-based models. The partial data-based models can use the sample
information sufficiently. Hereafter, we refer to this method as sufficient-sample-information
(SSI) method, which will be described in details in section 2.

In this paper, we propose weight choice criterion based on the Kullback-Leibler
(KL) loss. We prove that the average estimator is asymptotically optimal
in the sense that the corresponding KL loss is asymptotically identical
to that of the infeasible best model average estimator. To the best
of our knowledge, the current paper is the first work using model
averaging to deal with missing data in GLMs.

The rest of the article is organized as follows. Section 2 describes
the model framework and establishes the model averaging process. Section
3 shows the weight selection criteria and methods. Section 4 presents
asymptotic optimality. Section 5 reports the results of the simulation
studies on logistic and Poisson regression, respectively. Section
6 Summary. The relevant assumptions which are needed for the theorem's
proof and the proof processes are in the appendix.

\section{Model framework and model averaging process}

We consider the data generating process(DGP) 
\begin{equation}
f(y_{i}|\theta_{i},\phi)=exp\{\frac{y_{i}\theta_{i}-b(\theta_{i})}{\phi}+c(y_{i},\phi)\},i=1,\ldots,n.\label{1}
\end{equation}
where $\theta_{i}$ and $\phi$ are parameters, and $b(\cdot)$ and
$c(\cdot,\cdot)$ are known functions. The canonical parameter $\theta_{i}$
connects the parameter $\beta$ and the $K$-dimension covariate vector
$x_{i}$ in the form $\theta_{i}=x_{i}^{T}\beta$. Here we assume
that $K$ is fixed. Suppose we have $S$ candidate models and $S$
is finite. We estimate $\beta$ under different candidate models by
maximum likelihood estimation. Let $\theta_{0,i}$ be the true value
of $\theta_{i}$. We do not require that the true value $\theta_{0,i}$
is indeed a linear combination of $x_{i}$. In other words, it is
not required that there exist a $\beta$ so that $\theta_{0,i}=x_{i}^{T}\beta$.
Therefore, each of candidate model can be misspecified.

For a data set with $K$ covariates and the sample size $n$ to be
fixed, the number of possible covariate missing cases is $2^{K}-1$
(we collect covariates whose position of missing are same into one
group). Noting that not all such possible cases need be present in
a data set, we assume that we have $S-1\leq2^{K}-1$ covariate missing
cases, indexed by $s=2,\ldots,S$. For clear illustration, we provide
the following example.
\begin{description}
\item [{Example}]~
\end{description}
Suppose we have a sample with $K$ covariates and the sample size
is $n$. The covariates are missing completely at random. By sorting
the position of missing covariates, the following data matrix can
be obtained after a certain rearrangement of rows as well as columns:

\[
X=\left(\begin{array}{ccccc}
X_{11} & X_{12} & X_{13} & X_{14} & X_{15}\\
X_{21} & missing & missing & X_{24} & X_{25}\\
X_{31} & X_{32} & X_{33} & missing & X_{35}\\
X_{41} & X_{42} & missing & missing & missing
\end{array}\right)_{n\times K}
\]
where $X_{i_{1}i_{2}}(i_{1}\in\{1,2,3,4\},i_{2}\in\{1,2,3,4,5\})$
can be matrices. There are 5 cases with incomplete covariates: 
\begin{align*}
(X_{11}\ X_{21}\ X_{31}\ X_{41})^{T},\\
(X_{12}\ missing\ X_{32}\ X_{42})^{T},\\
(X_{13}\ missing\ X_{33}\ missing)^{T},\\
(X_{14}\ X_{24}\ missing\ missing)^{T},\\
(X_{15}\ X_{25}\ X_{35}\ missing)^{T},
\end{align*}
,so $S=6$ in this example.

We denote the CC model by model 1, under which all the covariates
are utilized but the sample size is smaller than $n$. Alternatively,
to utilize sufficient sample size, we may ignore some missing values
of covariates. Since we have $S-1$ cases of missing covariates, we
can have $S$ models in all, including one CC model and $S-1$ SSI
models. For the $s^{th}$ model, let $y^{(s)}$ be the associated
dependent variable and $X^{(s)}$ be the $n_{s}\times K_{s}$ covariate
matrix, where $n_{s}$ is the sample size, $K_{s}$ is the number
of covariate. We assume $X^{(s)}$ to be of full column rank. In example
1, model 1 has covariate matrix $(X_{11}\ X_{12}\ X_{13}\ X_{14}\ X_{15})$.
The other five SSI models have following covariate matrices: 
\begin{align*}
X^{(1)} & =(X_{11}\ X_{21}\ X_{31}\ X_{41})^{T},\ X^{(2)}=(X_{12}\ X_{32}\ X_{42})^{T},\\
X^{(3)} & =(X_{13}\ X_{33})^{T},\ X^{(4)}=(X_{14}\ X_{24})^{T}\ and\ X^{(5)}=(X_{15}\ X_{25}\ X_{35})^{T}.
\end{align*}

Let $\zeta_{s}$ denote the index set of the columns of X used in
model $s$. Let $\pi_{s}$ denote the projection matrix mapping $\beta=(\beta_{1},\ldots,\beta_{K})^{T}$
to the subvector $\pi_{s}\beta=\beta^{(s)}$ of components $\beta_{k},k\in\zeta_{s}$.
Denote the maximum likelihood estimator of $\beta^{(s)}$ as $\tilde{\beta}^{(s)}$.
Then the estimator of $\beta$ is ${\hat{\beta}}^{(s)}={\pi_{s}}^{T}{\tilde{\beta}}^{(s)}$
for the $s^{th}$ model. Some components of ${\hat{\beta}}^{(s)}$
are zeros.

Let $w=(w_{1},\ldots,w_{S})^{T}$ belonging in set: $\mathscr{H}=\{w\in[0,1]^{S}:\sum\limits _{s=1}^{S}w_{s}=1\}$.
Then the model averaging estimator of $\beta$ is: 
\[
\hat{\beta}(w)=\sum\limits _{s=1}^{S}w_{s}{\hat{\beta}}^{(s)}.
\]
Replace missing values in X by zeros and denote the resulting matrix
by $\tilde{X}$. In Example 1: 
\[
\tilde{X}=\left(\begin{array}{ccccc}
X_{11} & X_{12} & X_{13} & X_{14} & X_{15}\\
X_{21} & 0 & 0 & X_{24} & X_{25}\\
X_{31} & X_{32} & X_{33} & 0 & X_{35}\\
X_{41} & X_{42} & 0 & 0 & 0
\end{array}\right)_{n\times K}
\]
Suppose $\tilde{X}$ to be of full column rank, $y=(y_{1},\ldots,y_{n})^{T}$,
$\theta=(\theta_{1},\ldots,\theta_{n})^{T}$. $\theta_{0}$ is the
true value of $\theta$. Then a model averaging estimator of $\theta_{0}$
is: 
\[
\theta\{\hat{\beta}(w)\}=[\theta_{1}\{\hat{\beta}(w)\},\ldots,\theta_{n}\{\hat{\beta}(w)\}]^{T}=\tilde{X}\hat{\beta}(w).
\]
Note that $\tilde{X}$, not X, appears in $\theta\{\hat{\beta}(w)\}$,
because the use of X is infeasible. 

\section{Weight Choice}

Our weight choice criterion is: 
\begin{eqnarray}
%\nonumbertoremovenumbering(beforeeachequation)
\mathscr{G}(w) & = & 2\phi^{-1}B\{\hat{\beta}(w)\}-2\phi^{-1}y^{T}\theta\{\hat{\beta}(w)\}+\lambda_{n}w^{T}k,
\end{eqnarray}
where $\lambda_{n}=2$ like that of the panalty term of AIC, $K=(k_{1},\ldots,k_{S})^{T},\ k_{s}$
is the number of columns of X used in the $s^{th}$ candidate model.

$\mathscr{G}(w)$ comes from the Kullback-Leibler(KL) loss which is
defined as follows. Let $\mu=Ey$, $B_{0}=\sum\limits _{i=1}^{n}b(\theta_{0i})$,
$B\{\hat{\beta}(w)\}=\sum\limits _{i=1}^{n}b[\theta_{i}\{\hat{\beta}(w)\}]$,
and 
\[
J(w)=\phi^{-1}B\{\hat{\beta}(w)\}-\phi^{-1}\mu^{T}\theta\{\hat{\beta}(w)\}.
\]
The KL loss of $\theta\{\hat{\beta}(w)\}$ is 
\begin{eqnarray}
%\nonumbertoremovenumbering(beforeeachequation)
KL(w) & = & 2\sum\limits _{i=1}^{n}E_{y^{*}}\{log\{f(y^{*}|\theta_{0},\phi)\}-log\{f(y^{*}|\theta\{\hat{\beta}(w)\},\phi)\}\}\nonumber \\
 & = & 2\phi^{-1}B\{\hat{\beta}(w)\}-2\phi^{-1}\mu^{T}\theta\{\hat{\beta}(w)\}-2\phi^{-1}B_{0}+2\phi^{-1}\mu^{T}\theta_{0}\nonumber \\
 & = & 2J(w)-2\phi^{-1}B_{0}+2\phi^{-1}\mu^{T}\theta_{0},
\end{eqnarray}
where $y^{*}$ is another realization from $f(\cdot|\theta_{0},\phi)$
and independent of y. Assume $\phi$ is known. Typically, in logistic
and Poisson regressions, $\phi=1$. If $\mu$ was known, we could
obtain a weight vector by minimizing $J(w)$ given the relationship
between $J(w)$ and $KL(w)$ in (3). In practice, the minimization
of $J(w)$ is infeasible owing to the unknown parameter $\mu$. An
intuitive solution is to estimate $J(w)$. That is, we may use y to
estimate $\mu$ directly, i.e., we plug y into $J(w)$. Then we can
obtain weights by minimizing $\phi^{-1}B\{\hat{\beta}(w)\}-\phi^{-1}y^{T}\theta\{\hat{\beta}(w)\}$.
Unfortunately, this intuitive procedure leads to overfitting. To avoid
the overfitting, we use (2) as our weight choice criterion.

The resultant weight vector is defined as 
\begin{eqnarray}
%\nonumbertoremovenumbering(beforeeachequation)
\hat{w} & = & {argmin}_{w\in\mathscr{W}}\mathscr{G}(w)
\end{eqnarray}

\section{Asymptotic Optimality}

Let $\beta_{(s)}^{*}$ be the parameter vector which minimizes the
KL divergence between the true model (1) and the $s^{th}$ candidate
model. From Theorem 3.2 of \cite{white1982maximum}, we know that,
under certain regularity conditions, for $s\in\{1,\ldots,S\}$, 
\begin{eqnarray*}
%\nonumbertoremovenumbering(beforeeachequation)
\hat{\beta}_{(s)}-\beta_{(s)}^{*} & = & O_{p}(n_{s}^{-1/2}).
\end{eqnarray*}
Furthermore, if $\frac{n_{1}}{n}\textgreater c$, where $n_{1}$ is
the sample size of CC model and $c$ is a positive constant, since
$n>n_{s}>n_{1}$, we have $n_{s}=c^{*}n$ for some positive constant,
then
\begin{eqnarray}
%\nonumbertoremovenumbering(beforeeachequation)
\hat{\beta}_{(s)}-\beta_{(s)}^{*} & = & O_{p}(n^{-1/2}),
\end{eqnarray}
for $s\in\{1,\ldots,S\}$.

In order tostudy the optimality of the model averaging estimator,
we need the following conditions. 
\begin{description}
\item [{Condition(C.1)}] $||\tilde{X}^{T}\mu||=O(n),\ ||\tilde{X}^{T}\epsilon||=O_{p}(n^{1/2})$,
and uniformly for $w\in\mathscr{W}$, 
\begin{eqnarray*}
%\nonumbertoremovenumbering(beforeeachequation)
||\partial{B(\beta)}/\partial{\beta^{T}}|_{\beta=\tilde{\beta}(w)}|| & = & O_{p}(n)
\end{eqnarray*}
for every $\tilde{\beta}(w)$ between $\hat{\beta}(w)$ and $\beta^{*}(w)$.
\item [{Condition(C.2)}] Uniformly for $s\in\{1,\ldots,S\}$, $n^{-1}\overline{\sigma}^{2}||\theta(\beta_{(s)}^{*})||^{2}=O(1)$.
\item [{Condition(C.3)}] $n\xi_{n}^{-2}=o(1).$
\end{description}
The following theorem establishes the asymptotic optimality of the
model averaging estimator $\theta\{\hat{\beta}(w)\}$. \newtheorem{theorem}{Theorem}
\begin{theorem} If equation (5) and Conditions (C.1)-(C.3) are satisfied,
and $n^{-1/2}\lambda_{n}=O(1)$, then 
\[
{\frac{KL(\hat{w})}{{inf}_{w\in\mathscr{W}}KL(w)}\rightarrow1}\eqno{(6)}
\]
in probability as $n\rightarrow\infty$. \end{theorem} \newtheorem{proof}{Proof}
\begin{proof} See appendix. \end{proof}

Theorem 1 tell us, based on the estimated weight $\hat{w}$, the model
averaging estimate achieves the infimum of the KL loss.

\section{Simulation study}

In this section, we conduct two simulation experiments: logistic regression
and Poisson regression, to demonstrate the finite sample performance
of our model averaging method. Because our method can achieve asymptotic
optimality for missing data, we denote it by MOPT. In the simulations,
we compare the MOPT method with the CC model, model after mean imputation
(MIM), model averaging after mean imputation(MIMA). We set sample
size $n\in\{100,200\}$ and use KL loss (divided by n) for assessment.
For each setting, we generate 1000 simulated data. To mimic the situation
that all candidate models are misspecified, we pretend that the last
covariate missed in all candidate models. 

The simulation design is a logistic model. $y_{i}$ is generated from
$Binomial(1,p_{0i})$ with 
\[
p_{0i}=exp(x_{i}^{T}\beta)/\{1+exp(x_{i}^{T}\beta)\},
\]
where $\beta=(1,0.2,-1.2,-1,0.1)^{T}$ and $x_{i}=\left(x_{1i},x_{2i},x_{3i},x_{4i},x_{5i}\right)$
follow normal distribution with mean zeros, variance ones and the
correlations between different components of $x_{i}$ being 0.75.
In order to simulate covariates which are missing completely at random,
following \cite{ZHANG2013360}, we construct the missing data matrix
as follows. Data missing only occurs with $x_{i3}$ and $x_{i4}$
for some $i$. In order to control the missing structures, we generate
$\epsilon_{i}=(\epsilon_{i1},\epsilon_{i2})\sim N(0,I_{2})$, which
is independent with $x_{i}$. When $\epsilon_{ik}\textless a$, $x_{i(k+2)}$
is missed ($k=1,2$). We set the parameter $a\in\{-0.3,0,0.5\}$ to
control the ratio of missing observations. For MIMA, we consider all
possible submodels. 

Simulation results are as follows$(\times10^{-1})$:
\begin{center}
\begin{table}[H] \begin{center} \caption{Figure 1. Binomial}
\begin{tabular}{c c|c|c c c c}   
\hline\hline   
& n   &        & MOPT        & CC          & MIM       & MIMA   \\\hline   a=-0.3& 100 & mean   & 0.887       & 1.666       & 1.026     & 0.998    \\         &     & median & 0.867       & 1.169       & 0.899     & 0.959   \\         &     & SD     & 5.880       & 209.835     & 16.460    & 8.106  \\         & 200 & mean   & 0.732       & 0.635       & 0.774     & 0.800   \\         &     & median & 0.739       & 0.501       & 0.722     & 0.795    \\         &     & SD     & 3.114       & 14.688      & 0.549     & 3.834   \\\hline    a=0  & 100 & mean   & 1.068       & 13.907      & 1.369     & 1.231    \\         &     & median & 0.994       & 1.781       & 1.162     & 1.092    \\         &     & SD     & 1.755       & 72302.390   & 7.534     & 3.016    \\         & 200 & mean   & 0.906       & 1.091       & 1.017     & 0.995    \\         &     & median & 0.892       & 0.797       & 0.940     & 0.974   \\         &     & SD     & 3.700       & 65.761      & 7.864     & 4.420   \\\hline   a=0.5 & 100 & mean   & 1.321       & 274.463     & 1.775     & 1.450   \\         &     & median & 1.244       & 66.319      & 1.576     & 1.449  \\         &     & SD     & 0.532       &470643.900   & 1.791     & 0.484   \\         & 200 & mean   & 1.256       &15.346       & 1.549     & 1.380  \\         &     & median & 1.226       & 3.040       & 1.430     & 1.355  \\         &     & SD     & 0.329       &12859.850    & 1.471     & 0.573  \\   \hline \end{tabular} \end{center}\end{table}
\par\end{center}

Table 1. shows that when sample size $n$ increases or the proportion
of missing observations decreases, the mean and median values of the
KL loss of the MOPT decrease. The mean and median values of the other
four methods also decrease as the sample size $n$ increases or the
proportion of the missing observations decreases.

Except for the case $a=-0.3,n=200$, the mean and median values of
MOPT are both lower than all others. This pattern is also almost true
regarding standard deviation(SD) values except two cases of $a=-0.3,n=200$,
in which MOPT yields larger SD than MIM, and $a=0.5,n=100$, in which
MOPT yields larger SD than MIMA.

Moreover, the standard SD values of CC model are very huge comparing
with other models in most cases. This phenomenon show that the KL
loss is not stable for each replication of simulation for CC model.
This is becasue that each replication may have different data missing
patern that leads to significant differenc in sample size for the
CC model.

\section{Concluding Remarks}

In this paper, we propose a model averaging method (MOPT) to combine
the generalized linear models for the case when covariates are missing
completely at random to utilize the largest set of available cases.
The asymptotic optimality of our method has been proved. Developing
model averaging method to the generalized addictive models (GAMs)
with missing data and the nonlinear models is also worthy of study
in the future.

\section{Appendix}

All the limiting properties here and throughout the text hold under
$n\rightarrow\infty$. Let $\epsilon=(\epsilon_{1},\ldots,\epsilon_{n})^{T}=y-\mu$,
$\overline{\sigma}^{2}={max}_{i\in\{1,\ldots,n\}}var(\epsilon_{i})$,
$\beta^{*}(w)=\sum\limits _{s=1}^{S}w_{s}{\beta}_{(s)}^{*}$, 
\begin{eqnarray*}
%\nonumbertoremovenumbering(beforeeachequation)
KL^{*}(w) & = & 2\phi^{-1}B\{{\beta^{*}(w)}\}-2\phi^{-1}B_{0}-2\phi^{-1}\mu^{T}[\theta\{\beta^{*}(w)\}-\theta_{0}],
\end{eqnarray*}
and $\xi_{n}={inf}_{w\in\mathscr{W}}KL^{*}(w)$.\\

\[
Proof\ of\ Theorem\ 1
\]
Let $\tilde{\mathscr{G}}(w)=\mathscr{G}(w)-2\phi^{-1}B_{0}+2\phi^{-1}\mu^{T}\theta_{0}$.
Obviously, $\hat{w}={argmin}_{w\in\mathscr{W}}\tilde{\mathscr{G}}(w)$.
According to the proof of ${Theorem\ 1}'$ in Wan et al.(2010), Theorem
1 is valid if the following hold: 
\begin{flalign}
 & \sup\limits _{w\in\mathscr{W}}\frac{|KL(w)-KL^{*}(w)|}{KL^{*}(w)}=o_{p}(1) & \tag{A.1}
\end{flalign}
and 
\begin{flalign}
 & \sup\limits _{w\in\mathscr{W}}\frac{|\tilde{\mathscr{G}}(w)-KL^{*}(w)|}{KL^{*}(w)}=o_{p}(1) & \tag{A.2}
\end{flalign}

By (5), we know that uniformly for $w\in\mathscr{W}$, 
\begin{flalign}
 & \hat{\beta}(w)-\beta^{*}(w)=\sum\limits _{s=1}^{S}w_{s}(\hat{\beta}_{(s)}-\hat{\beta}_{(s)}^{*})=O_{p}(n^{-1/2}) & \tag{A.3}
\end{flalign}
It follows from (A.3), Condition (C.1) and Taylor expansion that uniformly
for $w\in\mathscr{W}$, \begin{flalign} & |B\{\hat{\beta}(w)\}-B\{\beta^*(w)\}| \leq  ||\frac{\partial{B(\beta)}}{\partial{\beta^T}}|_{\beta=\tilde{\beta}(w)}|| ||\hat{\beta}(w)-\beta^*(w)||=O_p(n^{1/2}),& \notag \\ &\mu^T[\theta\{\hat{\beta}(w)\}-\theta\{\beta^*(w)\}] \leq||\mu^TX||||\hat{\beta}(w)-\beta^*(w)||=O_p(n^{1/2}),&\notag \end{flalign} and \begin{flalign} &\epsilon^T[\theta\{\hat{\beta}(w)\}-\theta\{\beta^*(w)\}]\leq||\epsilon^TX||||\hat{\beta}(w)-\beta^*(w)||=O_p(1),&\notag \end{flalign}where
$\tilde{\beta}(w)$ is a vector between $\hat{\beta}(w)$ and $\beta^{*}(w)$.

In addition, using the central limit theorem and Condition (C.2),
we know that uniformly for $w\in\mathscr{W}$, 
\begin{flalign}
 & \epsilon^{T}\theta\{\beta^{*}(w)\}=\sum\limits _{s=1}^{S}w_{s}\epsilon^{T}\theta(\beta_{(s)}^{*})=O_{p}{(n^{1/2})}. & \tag{A.4}
\end{flalign}
These arguments indicate that 
\begin{flalign}
 & \sup\limits _{w\in\mathscr{W}}|KL(w)-KL^{*}(w)|\nonumber \\
 & =\sup\limits _{w\in\mathscr{W}}|2\phi^{-1}(B\{\hat{\beta}(w)\}-B\{\beta^{*}(w)\})-2\phi^{-1}\mu^{T}(\theta\{\hat{\beta}(w)\}-\theta\{\beta^{*}(w)\})|\nonumber \\
 & \leq2\phi^{-1}\sup\limits _{w\in\mathscr{W}}|B\{\hat{\beta}(w)\}-B\{\beta^{*}(w)|+2\phi^{-1}\sup\limits _{w\in\mathscr{W}}|\mu^{T}(\theta\{\hat{\beta}(w)\}-\theta\{\beta^{*}(w)\})|\nonumber \\
 & =O_{p}(n^{1/2}) & \tag{A.5}
\end{flalign}
and 
\begin{flalign}
 & \sup\limits _{w\in\mathscr{W}}|\tilde{\mathscr{G}}(w)-KL_{*}(w)|\nonumber \\
 & =\sup\limits _{w\in\mathscr{W}}|2\phi^{-1}(B\{\hat{\beta}(w)\}-B\{\beta^{*}(w)\})\nonumber \\
 & \quad-2\phi^{-1}(y^{T}\theta\{\hat{\beta}(w)\}-\mu^{T}\theta\{\beta^{*}(w)\})+\lambda_{n}w^{T}k|\nonumber \\
 & \leq2\phi^{-1}\sup\limits _{w\in\mathscr{W}}|B\{\hat{\beta}(w)\}-B\{\beta^{*}(w)\}|\nonumber \\
 & \quad+2\phi^{-1}\sup\limits _{w\in\mathscr{W}}|y^{T}\theta\{\hat{\beta}(w)\}-\mu^{T}\theta\{\beta^{*}(w)\}|+\lambda_{n}w^{T}k\nonumber \\
 & \leq2\phi^{-1}\sup\limits _{w\in\mathscr{W}}|B\{\hat{\beta}(w)\}-B\{\beta^{*}(w)\}|\nonumber \\
 & \quad+2\phi^{-1}\sup\limits _{w\in\mathscr{W}}|\mu^{T}[\theta\{\hat{\beta}(w)\}-\theta\{\beta^{*}(w)\}]|+2\phi^{-1}\sup\limits _{w\in\mathscr{W}}|\epsilon^{T}\theta\{\beta^{*}(w)\}|\nonumber \\
 & \quad+2\phi^{-1}\sup\limits _{w\in\mathscr{W}}|\epsilon^{T}[\theta\{\hat{\beta}(w)\}-\theta\{\beta^{*}(w)\}]|+\lambda_{n}w^{T}k\nonumber \\
 & =O_{p}(n^{1/2})+\lambda_{n}w^{T}k & \tag{A.6}
\end{flalign}

From (A.5), (A.6), Condition(C.3) and $n^{-1/2}\lambda_{n}=O(1)$,
we can obtain (A.1) and (A.2). This complete the proof.

\end{document}